\documentclass[12pt]{article}
\usepackage[utf8]{inputenc}
\usepackage{fullpage,amsmath,amsthm,amssymb,graphicx,authblk,url,comment}

\graphicspath{{figures/}}

\title{Deep Lagrangian connectivity in the global ocean inferred from Argo floats}

\author[1]{Ryan Abernathey}
\author[2]{Christopher Bladwell}
\author[2]{Gary Froyland}
\author[2,3]{Konstantinos Sakellariou}
\affil[1]{Lamont Doherty Earth Observatory of Columbia University,
New York, USA}
\affil[2]{School of Mathematics and Statistics, University of New South Wales, Sydney NSW 2052, Australia}
\affil[3]{Nodes \& Links Ltd, Leof. Athalassas 176, Strovolos, Nicosia, 2025, Cyprus}

\date{\today}

\begin{document}

\maketitle

\begin{abstract}
We describe the application of a new technique from nonlinear dynamical systems to infer the Lagrangian connectivity of the deep global ocean.
We approximate the dynamic Laplacian using Argo trajectories from January 2011 to January 2017 and extract the eight dominant coherent (or dynamically self-connected) regions at 1500m depth.
Our approach overcomes issues such as sparsity of observed data, and floats continually leaving and entering the dataset;  only 10\% of floats record for the full six years.
The identified coherent regions maximally trap water within them over the six-year time frame, providing a distinct analysis of the deep global ocean, and relevant information for planning future float deployment.
While our study is concerned with ocean circulation at a multi-year, global scale, the dynamic Laplacian approach may be applied at any temporal or spatial scale to identify coherent structures in ocean flow from positional time series information arising from observations or models.
\end{abstract}

\section{Introduction}

The Argo program has transformed oceanography by vastly increasing the spatiotemporal coverage of ocean observations \cite{Roemmich2009,RiserEtAl2016,Jayne2017}.
The most widely used aspect of Argo data are the temperature / salinity (T/S) profiles, which can provide, for example, estimates of ocean heat content e.g.~\cite{RoemmichEtAl2015,RiserEtAl2016,ChengEtAl2017} and constraints for ocean data assimilation systems e.g.~\cite{ForgetEtAl2015}.
However, Argo floats also provide a secondary type of data: their drift velocity at 1000 dbar over their 10-day cycle.
The aim of this paper is to use the Argo drift velocities to probe the structure of the mid-depth lateral circulation using new analytical tools from dynamical systems theory.

Compared to T/S profiles, Argo drift velocities have been relatively under-utilized.
Lebedev~\emph{et al.}~\cite{lebedev2007velocity} produced the first dataset of global drift velocities, and Katsumata and Yoshinari~\cite{KatsumataYoshinari2010} used this data to provide an esimate of the mean flow field at 1000 dbar.
Gray and Riser~\cite{GrayRiser2014} used individual Argo drift velocity vectors as a reference level to derive, together with the Argo geostrophic shear, a climatology of absolute geostrophic velocities over the upper 2000 dbar.
Argo drift velocities have also been used to study mesoscale processes.
Chiswell~\cite{Chiswell2013} used these velocities to make estimates of eddy diffusivity based on the deep Eulerian velocity decorrelation timescale.

Roach~\emph{et al.}~\cite{Roach2016} made an important conceptual breakthrough by considering long-term 1000 dbar {\em Lagrangian trajectories}; by analyzing simulated Argo floats in an ocean model, they concluded that the 10-day profiling cycle, in which the floats return to the surface to transmit their data, does not significantly disrupt the Lagrangian trajectory at 1000 dbar.
In other words, simulated floats with profiling cycles had very similar trajectories and dispersion characteristics to those that did not profile but just remained at 1000 dbar.
This means that Argo floats can be used not only to study 10-day displacement vectors but also trajectories over much longer timescales.
This property was further exploited by
Roach~\emph{et al.}~\cite{Roach2016}, and subsequent studies by Balwada~\emph{et al.}~\cite{BalwadaEtAl2016,BalwadaEtAl2021}, to calculate single and multi-particle dispersion statistics.
The notion that Argo trajectories may represent unbroken water mass trajectories over long time intervals opens the door to many tools from dynamical systems \cite{dijkstra2005} and Lagrangian \cite{wiggins2005,van2018} analysis.


  Over the last two decades a prominent subset of Lagrangian methods have aimed at identifying coherent sets and structures in fluids, with ocean flow as an important motivating testbed.
  These objects have  influence on for example:  the motion of surface drifters \cite{olascoaga2013}, the accumulation of plastic on the ocean surface \cite{froyland2014well}, and the risk of bycatch \cite{scales2018}.
   One prevalent class of Lagrangian methods uses linearisations of the flow to identify curves (e.g.\ on the surface ocean) that are distinguished in some way:  usually as local extremisers of local spatial expansion, possibly constrained tangentially or transversally to the curve.
   These include finite-time Lyapunov exponents (FTLEs) \cite{pierrehumbert1994tracer, haller2000lagrangian,waugh2012}, and variational theory \cite{hallervariational,haller2013coherent}, to name just a few of these approaches.
    The vast majority of these techniques rely on accurate spatial derivatives of the flow, which are impossible to obtain from Argo drifter observations because of the latter's sparsity.

 Instead of using spatial derivatives of the nonlinear flow, transfer operator approaches \cite{FSM10,F13} create a global linear description of the flow (the transfer operator) and extract coherent regions from singular vectors of this operator.
 Finite-time coherent sets were introduced in \cite{FSM10} as  time-parameterised families of regions in a fluid that have \emph{persistently small leakage} from the family over the finite-time evolution.
    The transfer operator method for estimating finite-time coherent sets \cite{FSM10} has been applied in the ocean setting to identify ocean eddies in three dimensions \cite{froylandetal12} and study their surface decay \cite{froylandetal15}.
    This approach directly estimates the coherent regions themselves, e.g.\ the eddy as a two-dimensional or three-dimensional object on the surface ocean or in the three-dimensional ocean (or in two-dimensional surfaces or three-dimensional volumes in the atmosphere \cite{FSM10}).
    Because the transfer operator method \cite{FSM10} does not require spatial derivatives, it is robust to noise, and has less stringent requirements on the resolution of trajectory sampling.

    In a move toward sparser data requirements, Froyland and Padberg-Gehle \cite{FPG15} showed fuzzy $c$-means clustering could extract coherent sets from sparse data sets, and provided a 5-year coherent set analysis of the surface ocean from the Global Drifter Array (GLAD), similar to the analysis we will describe in this work.
   A diffusion-maps implementation of the transfer operator method \cite{F13} was described in \cite{BK17}, who were able to reproduce the surface ocean results from \cite{FPG15} using the same trajectory dataset.

    In the present work we use the dynamic Laplacian \cite{F15} approach, and in particular the finite-element method (FEM) implementation \cite{FJ18}.
    The dynamic Laplacian formally arises as a ``zero-noise limit'' of the transfer operator approach \cite{F13}, and therefore one expects very similar estimates of finite-time coherent sets.
    In the next section, we discuss an equivalent characterisation of finite-time coherent sets as a time-parameterised family of sets whose \emph{boundary remains persistently small} (i.e.\ does not become filamented, and therefore minimises diffusive leakage) under evolution.
    Both the transfer operator and dynamic Laplacian approaches are a form of spectral clustering, adapted to the Lagrangian dynamics we are investigating.

    The FEM implementation \cite{FJ18}  handles sparse, scattered data well, including piecing together trajectory segments formed by short drifter trajectories that do not last for the full time interval of analysis (adapting an approach from \cite{FPG15}).
    The FEM approach also has a number of advantages over other implementations such as \cite{hadji16,BK17}, for instance (i) being unaffected by varying spatial densities of trajectories, (ii) providing full coherent set ``fields'' on the domain (not only at trajectory points), and (iii) not requiring any parameter tuning.
    Because of the above points, we use this FEM--dynamic Laplacian approach in this work to analyse the ocean at the $\sim$1000 metre level, using only the sparsely distributed Argo float position information encoded in their trajectories.

    Related analysis of \emph{almost-invariant sets}, namely regions of the ocean that remain approximately \emph{fixed in space} under the ocean flow, have been undertaken using float trajectories in the surface Gulf of Mexico \cite{gomsurface} using a combination of drifters from several drifter arrays, and at 1500--2500m depth in the Gulf of Mexico \cite{gomdeep} using RAFOS floats.
    In these studies, the underlying ocean currents were assumed to be stationary in a statistical sense, meaning that the likelihood of drifter movements between similar regions of the ocean at different times did not change over time.
    This enabled a huge multiplicative ``reuse'' of the data (similar to \cite{vsef12}), and allowed the more classical Ulam approach for finding almost-invariant sets, e.g. \cite{FPET07} to be employed.

    As far as we are aware the current study is the first truly Lagrangian study of coherent sets that directly uses deep ocean float trajectories.
   The FEM implementation of the dynamic Laplacian \cite{FJ18} approach we use here is ideally suited to this type of transport analysis when only sparse, scattered trajectory data is available.
   This FEM approach has been previously used in the oceanographic context to track eddy motion using satellite-derived velocity fields \cite{FJ18,SEBA}, as well as to track the meander of the Gulf Stream from purely kinematic trajectory positional information \cite{SEBA};  these experiments were at smaller spatial scales and used richer, synthetic trajectory information compared to the Argo trajectories used here.

Our goals in this paper are twofold.
First, we wish to prototype the dynamic Laplacian methodology using real-world data and demonstrate that it is computationally practical and gives reasonable results.
We hope that this demonstration will inspire others to adapt this method to different datasets.
Second, we wish to use the results of the analysis to investigate the dynamic geography of the mid-depth ocean flow, a topic with multiple applications.
The coherent sets identified by this method have significant implications for the Argo deployment strategy, since they effectively reveal the regions that are poorly dynamically connected by Argo floats.
We can also explore the relationships between the Argo coherent sets and large-scale ocean tracers, although, as described below, these results are somewhat inconclusive.

An outline of the paper is as follows.
In section 2 we give a brief background on finite-time coherent sets and define the dynamic Laplace operator.
Section 3.1 describes Argo drifter data and the method by which we create trajectories from the drifter movements.
Section 3.2 details the numerical discretisation of the dynamic Laplacian and how we compute this discretisation using trajectories from Argo drifters.
In section 3.3 we describe how to isolate individual coherent sets from the eigenfunction of the dynamic Laplacian using the Sparse EigenBasis Approximation (SEBA) algorithm.
Section 4 opens with a discussion of the eight most coherent regions we identify in the context of the ocean circulation at 1000m depth.
In section 4.1 we describe how to find precise boundaries for the eight maximally coherent sets.
Section 5 compares the locations of the eight identified global-scale coherent sets with the large-scale properties planetary potential vorticity and oxygen concentration.
We conclude in section 6 with a summary of our findings, the advantages of the dynamic Laplacian analysis approach, and the relevance of coherent set analysis to elucidating ocean circulation properties and to informing the placement of drifters in programs such as Argo.



\section{Finite-time coherent set analysis and the dynamic Laplacian}

Using trajectories of the Argo float array, we wish to identify coherent regions in the deep global ocean.
The trajectories of the Argo floats at 1000m depth induce a nonautonomous dynamical system on a two-dimensional phase space at this depth.
Finite-time coherent sets were introduced in \cite{FSM10,F13} as subsets of phase space that \emph{minimally mix with (or leak into) the rest of the  phase space}.
Coherent features such as mesoscale eddies \cite{froylandetal12,froylandetal15} and the Antarctic polar vortex \cite{FSM10} have been identified using transfer operator methods using this least-mixing / least-leaking criterion.

An alternate, but strongly related criterion for finite-time coherence is based on geometric considerations \cite{F15}.
In the presence of small diffusive processes, the extent to which a region will mix with the surrounding phase space is proportional to average boundary size over the time interval of interest.
Let $X$ denote the 1000m depth level (with land removed);  this is a subset of $\mathbb{R}^2$.
By $[0,t_f]$ we denote the time interval over which we analyse the float dynamics, where we nominally begin at an initial time of 0.
Let $\phi:[0,t_f]\times X\to X$ denote the nonlinear time-dependent flow map that evolves a location $x\in X$ at time $0$ to its future location $\phi_t(x)$ at time $t$.
We note that in our application in this paper, the domain $X$ is invariant under the flow, i.e.\ $\phi_t(X)=X$, but this is not necessary for the methodology to be applied.
Let $A\subset X$ and we denote by $\partial A$ and $\partial X$ the boundaries of $A$ and $X$ respectively.
We will be concerned with regions $A$ at the 1000m depth level so that $A\cap\partial X=\emptyset$.
Let $\ell(\partial A)$ denote the length of the boundary of $A$ and $a(A)$ denote the area\footnote{If $X$ were three-dimensional, $\ell(\partial A)$ and $a(\partial A)$ would be replaced with the area and volume, respectively, of the boundary of $A$.} of $A$.
Finite-time coherent sets are sets $A\subset X$ that produce close to minimal values of
\footnote{The quantity $h(A)$ is the dynamic Cheeger value of $A$ \cite{F15}, written here in the ``Dirichlet boundary condition'' version;  see \cite{FJ18} for details.}
\begin{equation}
    \label{cheeger}
    h(A):=\frac{1}{t_f}\int_0^{t_f} \frac{\ell(\partial (\phi_t(A)))}{a(\phi_t(A))}\ dt.
\end{equation}
Thus, we seek $A\subset X$ whose time-averaged (under forward evolution of $A$ by $\phi$) ratio of boundary length to enclosed area is minimised.
Such a set will have minimal mixing in the presence of small diffusion.
\begin{figure}[hbt]
    \centering
    \includegraphics[width=.75\textwidth]{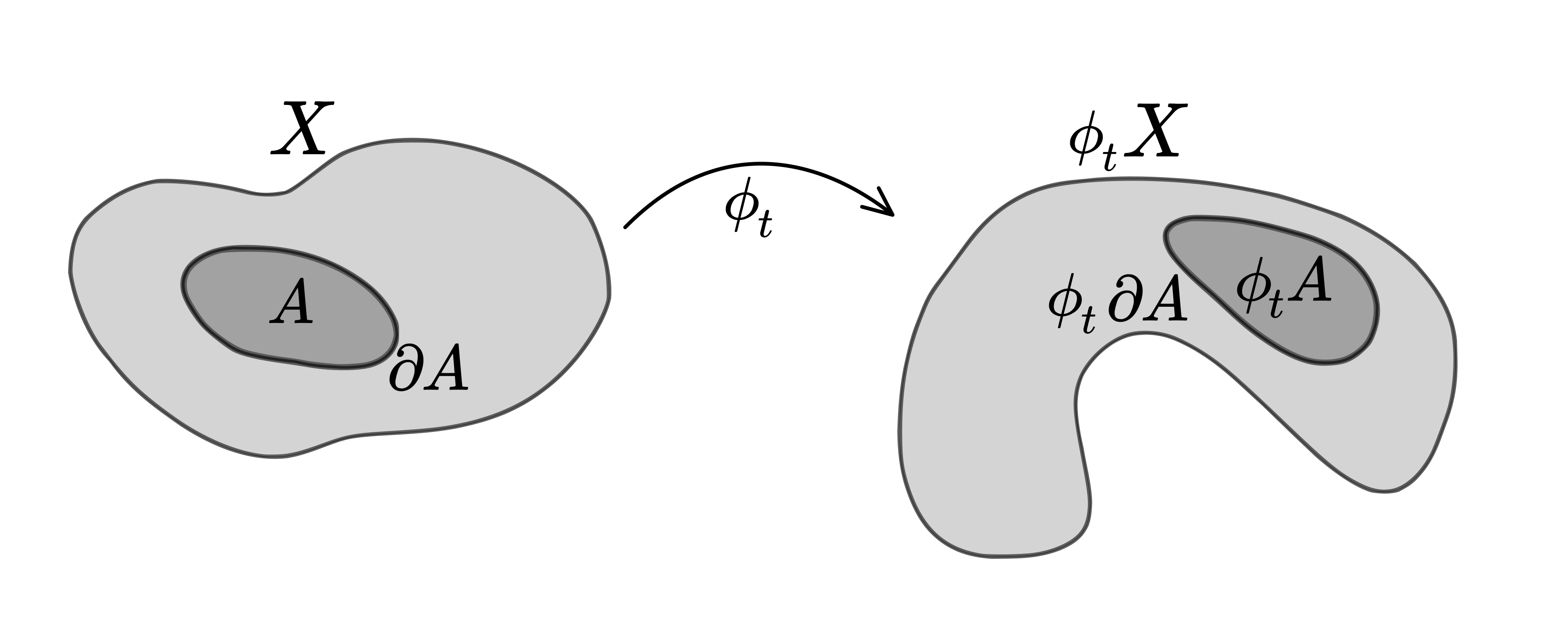}
    \caption{Illustration of a region of the ocean (light gray), containing a set $A$ (dark gray), whose boundary does not grow significantly under forward evolution.}
    \label{fig:CS}
\end{figure}

To find these regions, we will use the \emph{dynamic Laplace operator} \cite{F15}, which we now describe.
We denote by $\Delta$ the usual Laplace operator acting on functions $f:X\to \mathbb{R}$ defined on our 1000m depth level $X$.
It is well known that the spectrum and eigenfunctions of $\Delta$ encode geometric features of the domain $X$ \cite{SIAMreview}.  This is related to the question posed in the 1966 paper of Mark Kac:  ``Can one hear the shape of a drum?'' \cite{kac66}, which asks if a surface can be classified by the spectrum of the Laplace operator on that surface.
Further, the spectrum and eigenfunctions of $\Delta$ appear in {\emph{isoperimetric theory}} (see for example \cite{chaveleigen}), which concerns sets with minimal boundary size relative to volume.

The dynamic Laplace operator extends these ideas to dynamical systems to answer \emph{dynamic} isoperimetric questions like the minimal values of (\ref{cheeger}).
Let $f:X\to\mathbb{R}$ denote a function that we think of as being ``anchored'' at time $0$.
There is a natural way to push forward this function under the dynamics, namely by precomposing with the inverse flow map $\phi_{-t}$.
That is, the function $f\circ\phi^{-t}:X\to\mathbb{R}$ is anchored at a time $t$ in the future and strictly speaking is a function defined on $\phi^t(X)$, but in our deep global ocean application, we have $\phi^t(X)=X$ for all $t\in[0,t_f]$.
To access the geometry of $\phi^t(X)$ (at time $t$), which has been deformed by the nonlinear $\phi^t$, we apply the usual Laplace operator $\Delta$ to the function $f\circ \phi^{-t}$, which is anchored at time $t$.
Having formed $\Delta(f\circ\phi^{-t})$ (which is just another function anchored at time $t$), we need to pull it back to the initial time 0.
We do this by precomposing with $\phi^t$, leading finally to $(\Delta(f\circ\phi^{-t}))\circ\phi^t$ as the contribution to the dynamic Laplacian from time $t$.
We average these contributions across the interval $[0,t_f]$ to obtain
\begin{equation}
    \label{dynlap}
    \Delta^D_{[0,t_f]}:=\frac{1}{t_f}\int_0^{t_f} (\Delta(f\circ\phi^{-t}))\circ\phi^t\ dt=\frac{1}{t_f}\int_0^{t_f} (\phi_t)^*\circ\Delta\circ(\phi_t)_*\ dt,
\end{equation}
where the second equality uses the standard push-forward and pull-back notation $(\phi_t)_*f:=f\circ(\phi_t)^{-1}$ and $(\phi_t)^*f:=f\circ\phi_t$.
Just as the eigenfunctions of the standard Laplace operator $\Delta$ contain information about regions with minimal boundary length relative to enclosed area, the leading eigenfunctions of the dynamic Laplace operator $\Delta^D_{[0,t_f]}$ encode regions with \emph{minimal average boundary length under evolution with $\phi_t$}, namely \emph{finite-time coherent sets}.
Using Argo float trajectories, we will infer $\phi_t$, numerically construct $\Delta^D_{[0,t_f]}$, and compute its leading eigenfunctions.

\section{Approximating dynamic Laplacian eigenfunctions from Argo data}

\subsection{Description of the Argo-based trajectory data}
\begin{figure}[hbt]
    \centering
    \includegraphics[width=0.49\textwidth]{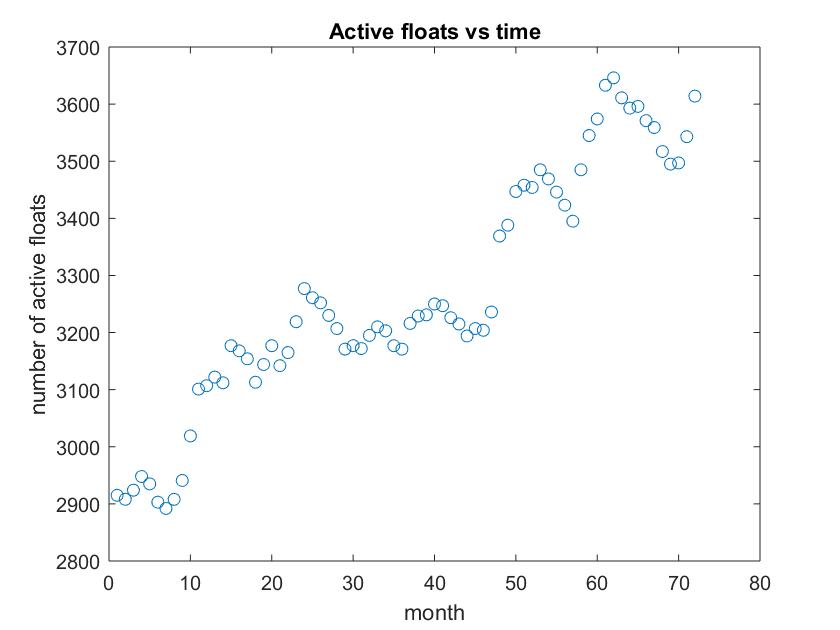}
    \includegraphics[width=0.49\textwidth]{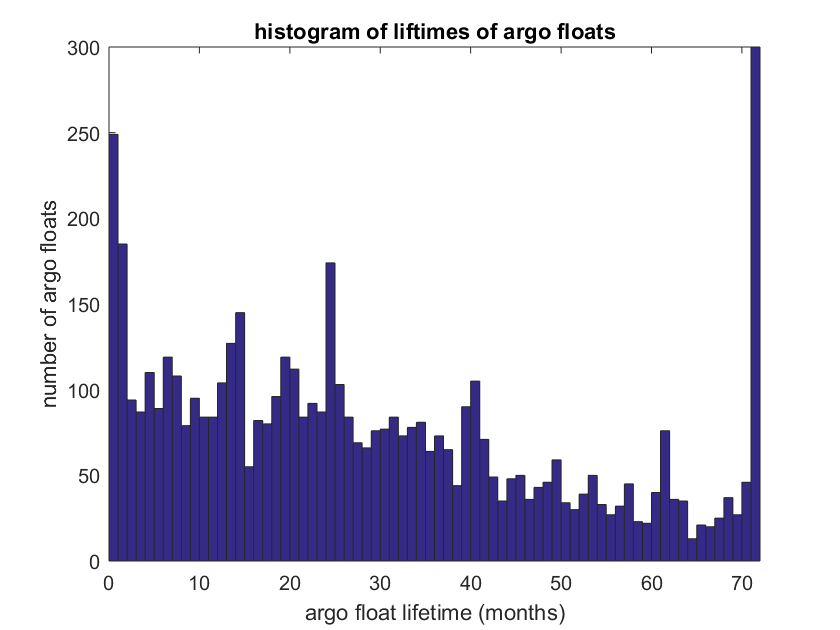}
    \caption{Left: Active floats vs.\ month.  Right: Histogram of float lifetimes over six years.}
    \label{fig:floatdata}
\end{figure}

The Argo program is an international collaboration which began in the late 1990s with the goal of providing freely available subsurface temperature and salinity measurements to understand ocean climate variability.
By deploying a fleet of autonomous devices called floats, the Argo program has been able to capture measurements in the upper 2000m of the water column over most of the global ocean.
Currently, Argo is the sole global subsurface data set.
The floats are approximately one metre in length and are able to alter their density by extruding mineral oil contained inside them into an external bladder.
This mechanism allows them to descend to a predetermined pressure, called the parking pressure, which allows them to drift at a depth of approximately 1000m below sea level.
Figure \ref{fig:floatdata} shows the number of Argo floats that were active in the 72 months from January 2011 to January 2017, and the distribution of lifetimes of the Argo floats.

Floats in the Argo program follow cycles that involve an extended period at a depth followed by a brief period at the surface.
A typical approximately 10-day cycle begins as the device descends from the surface to its parking pressure at around 1000m, where it spends on average 9 days.
The device then descends to a maximum depth of 2000m, after which it returns to the surface.
The average displacement of each of the floats during the subsurface component of each cycle is 37.64km and each float completes on average 140 cycles over its lifespan.
During the surface component of each cycle, lasting an average of 8 hours, the device transmits its location on average 9 times, with an accuracy of approximately 100m.
We will use the initial location of the float upon surfacing.
The raw netCDF data is freely available from the GDAC websites  (\url{http://www.coriolis.eu.org/} or \url{http://www.usgodae.org/Argo/}).
Full specifications can be found in \cite{Argo}.

The Argo floats are deployed individually, and their surfacing times are not coordinated.
We create a set of Lagrangian trajectories, approximately recorded monthly as follows.
There are $I=8140$ unique floats that are active at \emph{some} time in the 6-year period under consideration; we index these floats by $i\in\{1,2,\ldots,I\}$. 
For each month, indexed by $t\in\{1,2,\ldots,T\}$, $T=72$, we record the surfacing locations of floats that surface between the 1st and 12th days of that month.
We then create trajectories of monthly float locations $\{x_{i,t}\}_{1\le i\le I,1\le t\le T}$ by
\begin{equation}
    \label{trajdef}
x_{i,t}=\left\{
  \begin{array}{ll}
    \mbox{first surfacing location}, & \hbox{if float $i$ surfaces between the 1st and 12th of month $t$;} \\
    \emptyset, & \hbox{otherwise.}
  \end{array}
\right.
\end{equation}
For later use we define the set of reporting floats $R_t$ at time $t$ by $R_t=\{i:x_{i,t}\neq \emptyset\}$, $t=1,\ldots,T$.
We will additionally use a discrete set of $C=2712$ points $\{y_j\}_{1\le j\le C}$ sampled on the coastlines of all continents and large islands.
These $C$ points were obtained from the \textsc{Matlab} Mapping Toolbox's \verb"worldmap" dataset, subsampling every 5th point and manually ensuring large islands were well represented.
We emphasise that we use each space/time-stamp of each float only once in our computations, in contrast to e.g.\ \cite{vsef12,gomsurface,gomdeep} where a float position at one time is reused many times over.



\subsection{Numerical approximation of dynamic Laplacian and its eigenproblem}\label{section: FEM dyn lap}
Techniques of numerically approximating the dynamic Laplacian eigenproblem include \cite{F15,froyland2015fast,BK17};  related constructions include \cite{hadji16} and some aspects of \cite{padberg2017}.
To handle the Argo data we use the ``adaptive TO'' (adaptive transfer operator) approach from \cite{FJ18} based on the finite-element method.
In addition to having a very well-developed theoretical \cite{F15} and numerical \cite{FJ18} foundation, this approach has several practical advantages:
\begin{itemize}
    \item  the trajectories can be scattered in space and sparse (the Argo floats are sparsely distributed),
    \item the domain can have irregular boundaries (the ocean is bounded by irregular continents),
    \item complete trajectories are not required (Argo floats only report over some fraction of the timespan $[0,t_f]$),
    \item no spatial derivatives are required (this would be impossible with sparse Argo data), \item eigenfunctions are reconstructed on the full phase space, not only at data points,
    \item there are no parameters to select.
\end{itemize}


We consider the eigenproblem of the dynamic Laplacian with Dirichlet boundary conditions, namely $\Delta^D_{[0,t_f]} f=\lambda f$ on the interior of $X$ and $f\equiv 0$ on $\partial X$.
The dynamic Laplacian is a symmetric, elliptic operator \cite{F15,FJ18} and its eigenspectrum is real and nonpositive.
We are interested in the eigenfunctions $f$ corresponding to the large values of $\lambda$ (i.e.\ those negative eigenvalues closest to 0).

Using the central term of (\ref{dynlap}), this eigenproblem may be written as $$\frac{1}{t_f}\int_0^{t_f} (\Delta(f\circ\phi_{-t}))\circ\phi_t\ dt=\lambda f.$$
Multiplying through by a function $g:X\to\mathbb{R}$, integrating over $X$, and applying integration by parts to move one spatial derivative from $f$ to $g$, we have
the weak form of this eigenproblem:
\begin{equation}
    \label{weak1}
    -\frac{1}{t_f}\int_0^{t_f}\int_{X} \nabla(f\circ\phi_{-t})(x)\cdot\nabla(g\circ\phi_{-t})(x)\ dx\ dt = \lambda\int_{X} f(x)g(x)\ dx.
\end{equation}

\begin{figure}[hbt]
    \centering
    \includegraphics[width=\textwidth]{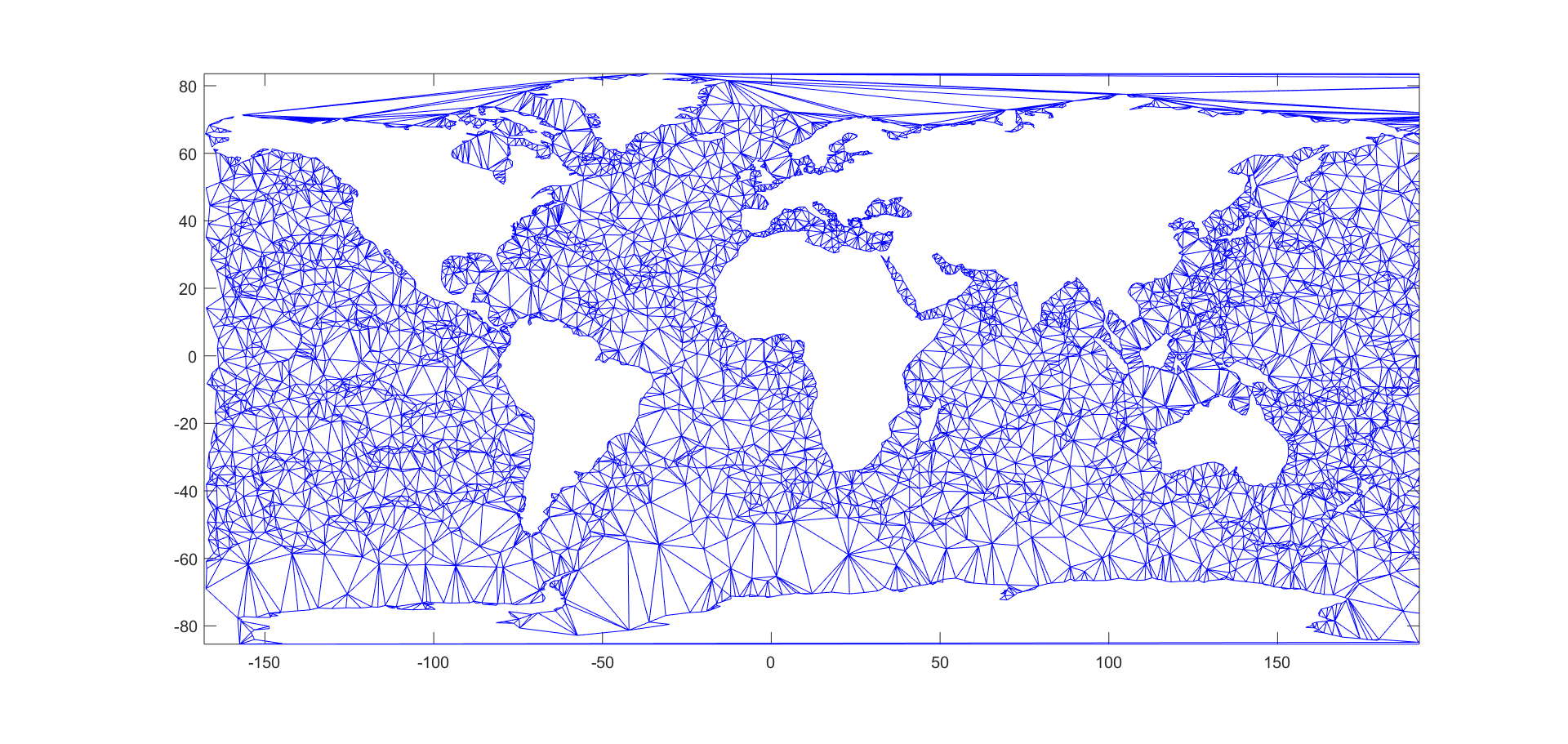}
    \caption{Mesh at $t=36$ months. Vertices of the mesh are either the positions of Argo floats at 36 months or points on coastlines.}
    \label{fig:mesh}
\end{figure}

\begin{figure}[hbt]
    \centering
    \includegraphics[width=0.5\textwidth]{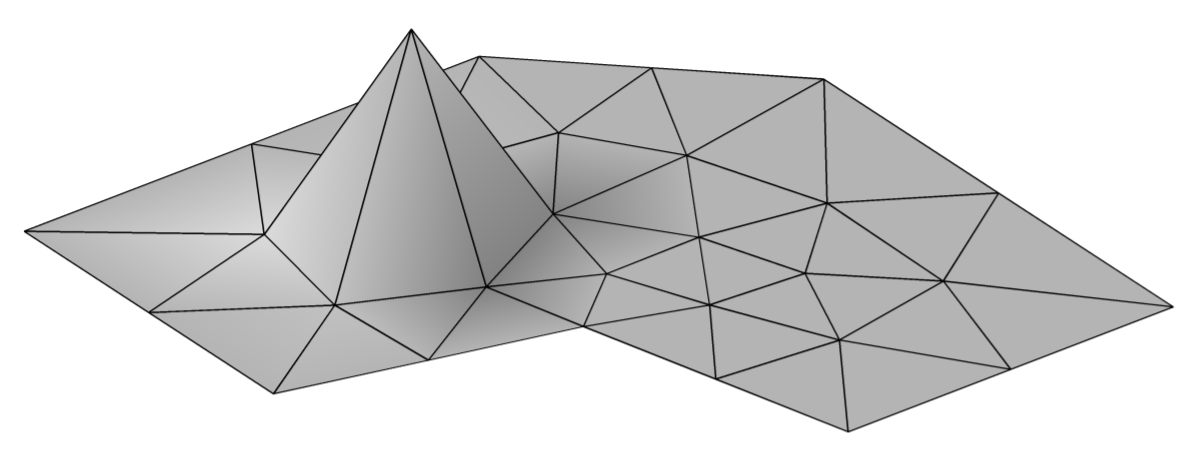}
    \caption{Example graph of a single piecewise linear hat function $\varphi^t_i$ on a two-dimensional mesh. The function $\varphi_{i}^t$ takes the value 1 at the single vertex $x_{i,t}$ and is zero at neighboring vertices (and all other vertices). Reproduced from \cite{arnold10}.}
    \label{hatfunction}
\end{figure}

Our numerical approximation will be based upon the float locations $\{x_{i,t}\}$ in (\ref{trajdef}) for a discrete set of $t$ in the interval $[0,t_f]$.
We follow Sections 3, 3.2, 3.2.3, and 3.2.4 in \cite{FJ18}.
The main preparatory steps are:
\begin{enumerate}
\item We discretise the interval $[t_0,t_f]$ into times indexed by $t=1,\ldots,T$.  In the computations presented here, the time index $t$ indicates an integer number of months, and $T=72$.
\item At each $t=1,\ldots,T$ we mesh the positions of the reporting floats and coastline points $\{x_{i,t}: i\in R_t\}\cup\{y_j: 1\le j\le C\}$; see Figure \ref{fig:mesh}. This mesh covers a region that is very close to $X$ and we will henceforth not distinguish between the domain covered by the mesh and $X$.
\item At each $t=1,\ldots,T$ using the mesh at time $t$, we build a collection (a linear basis) of piecewise linear hat functions $\{\varphi_1^t,\ldots,\varphi_{n_t}^t\}$,  where $n_t=|R_t|+C$.
We will approximate eigenfunctions of the dynamic Laplacian on $X$ by linear combinations of the $\{\varphi_i^1\}_{i=1}^{n_t}$.
A sketch of one of these hat functions is illustrated in Figure \ref{hatfunction}.
\end{enumerate}

By linearity of (\ref{weak1}) in $f$ and $g$ we may consider basis pairs $f=\varphi^t_i$ and $g=\varphi_j^t$, and for each $t=1,\ldots,T$ compute the $n_t\times n_t$ ``stiffness'' matrices
\begin{equation}
    \label{Dmatrix}
    D^t_{ij}=\int_{X}\nabla \varphi_i^t(x)\cdot\nabla\varphi_j^t(x)\ dx,
\end{equation}
for $i,j\in R_t$, setting $D^t_{ij}=0$ otherwise.
Similarly, for each $t=1,\ldots,T$ we compute ``mass'' matrices
\begin{equation}
    \label{Mmatrix}
    M^t_{ij}=\int_{X} \varphi_i^t(x)\varphi_j^t(x)\ dx,
\end{equation}
for $i,j\in R_t$, setting $M^t_{ij}=0$ otherwise.
We note that because the $\varphi_i^t$ are piecewise linear functions of space, the integrands in (\ref{Dmatrix}) and (\ref{Mmatrix}) are piecewise constant and quadratic, respectively and can be efficiently evaluated numerically \cite{FJ18}.

We now wish to sum the matrices $D^t$ and $M^t$ over $t$ and so for each $t$ we inflate each $D^t$ and $M^t$ from size $(|R^t|+C)\times (|R^t|+C)$ to size $(I+C)\times (I+C)$ by adding zero rows and columns at indices in $I\setminus R_t$.
For $u\in \mathbb{R}^{I+C}$, the discrete form of (\ref{weak1}) is the eigenproblem
\begin{equation}
    \label{discreteweak}
    -\left(\frac{1}{T}\sum_{t=1}^{T}D^t\right)u = \lambda\left(\frac{1}{T}\sum_{t=1}^{T}M^t\right)u,\qquad u_i=0\mbox{ for }i=I+1,\ldots,I+C,
\end{equation}
where we fix $u_i$ to 0 for all indices corresponding to coastline points in order to enforce Dirichlet boundary conditions (namely the last $C$ points according to our indexing).
In this study we have selected Dirichlet boundary conditions because we are seeking slowly mixing structures that do not intersect with coastlines.

\subsection{Eigenvector and finite-time coherent set computations}

\begin{figure}
  \centering
  \includegraphics[width=\textwidth]{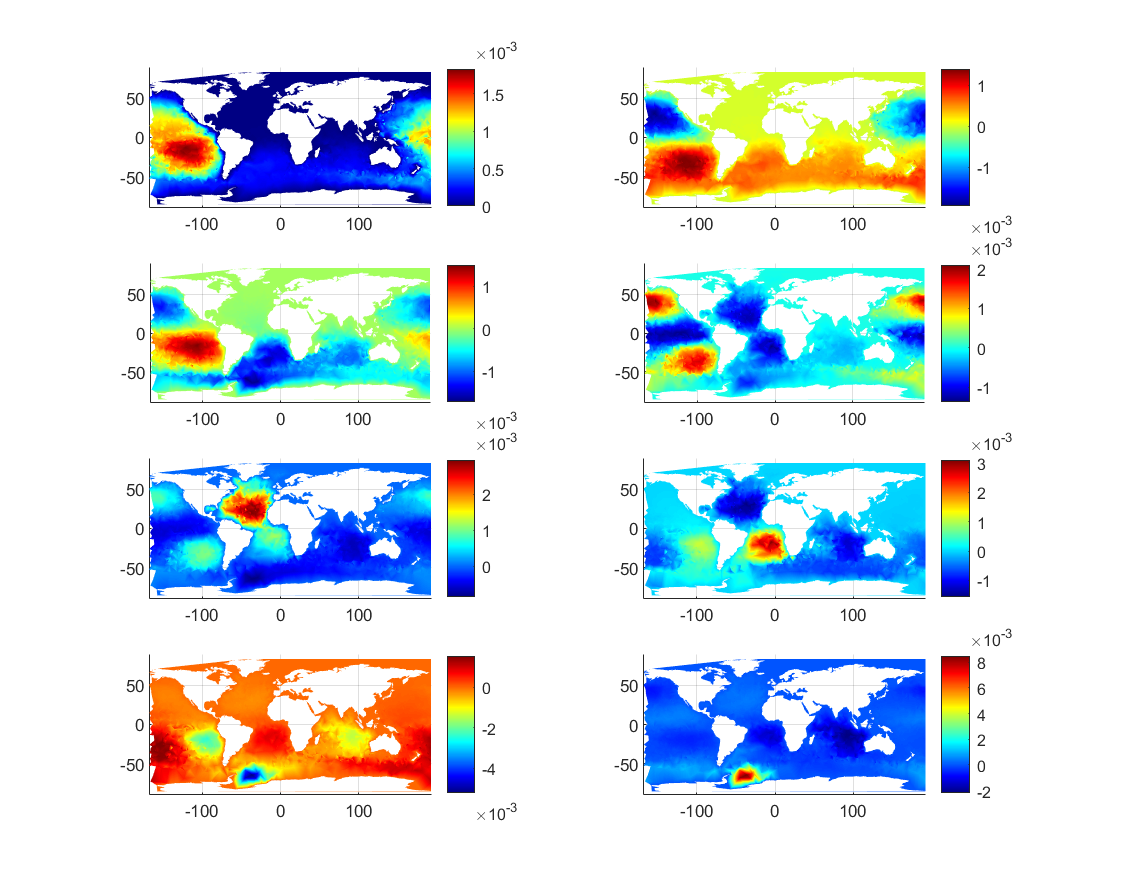}
\caption{Leading eight eigenfunctions $f^{(k)}$, $k=1,\ldots,8$, of the dynamic Laplacian computed from Argo trajectory data over six years. These eigenfunctions are represented as scalar fields at January 2011 (month \# 1 for our trajectories). Deep red and deep blue regions correspond to finite-time coherent sets.}
\label{fig:8evecs-hatfunc}
\end{figure}

We use (\ref{discreteweak}) to find eigenvalue, eigenvector pairs $(\lambda,u)$.
From the coefficient vector $u$, we approximate the corresponding eigenfunction $f$ by $f(x)=\sum_{i=1}^{I+C} u_i\varphi_i^1(x)$, which is anchored at time 0.
We compute the leading eight eigenvectors, denoted $u^{(k)}$, $k=1,\ldots,8$, and from these we construct the leading eight eigenfunctions $f^{(k)}=\sum_{i=1}^{I+C}u^{(k)}_i\varphi_i^1$;   see Figure \ref{fig:8evecs-hatfunc}.
We may estimate the forward evolved eigenfunction $f\circ\phi_{-t}$ at month $t$ by $(f\circ\phi_{-t})(x)\approx\sum_{i=1}^{I+C}u_i\varphi_i^t(x)$.

\begin{figure}
  \centering
  \includegraphics[width=\textwidth]{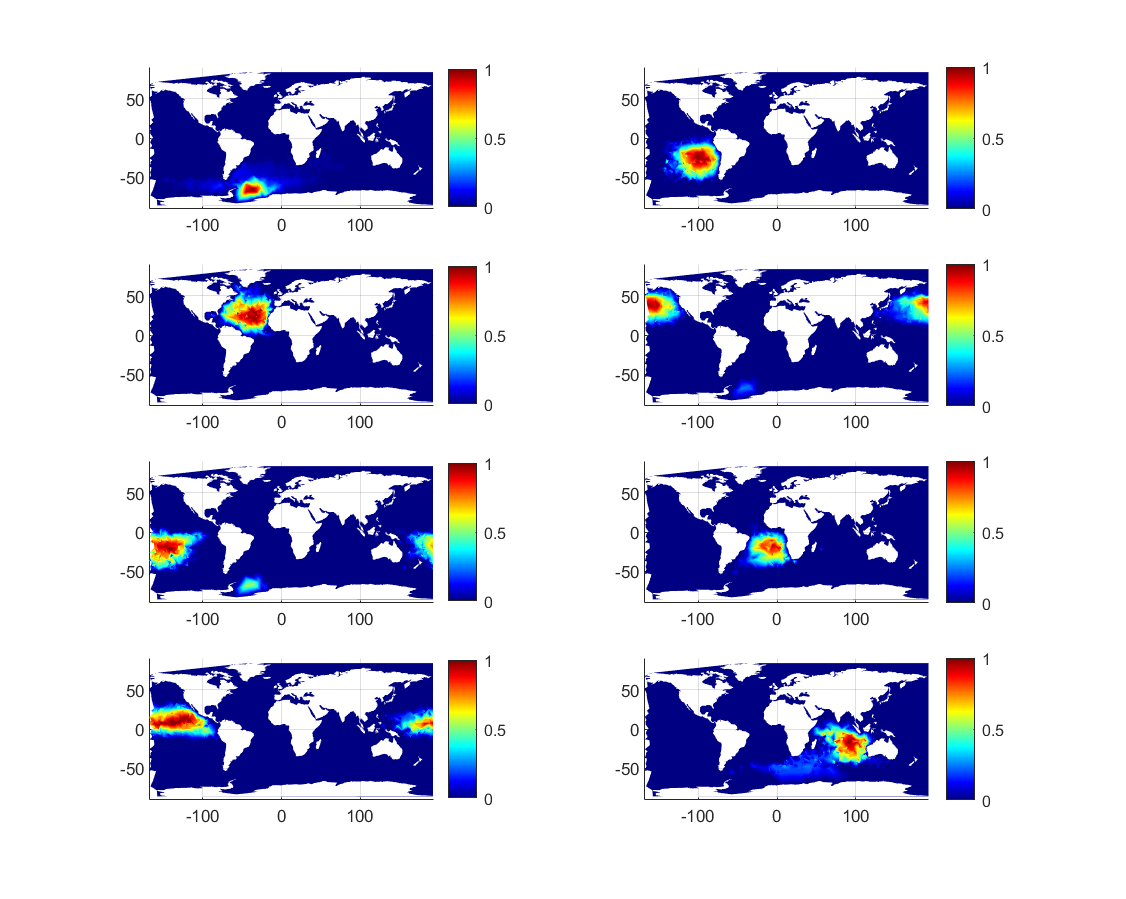}
\caption{Result of the application of the SEBA algorithm \cite{SEBA} to the leading 8 eigenvectors to separate the eight main coherent sets. These sparse basis functions $\hat{f}^{(k)}$, $k=1,\ldots,8$, are represented as scalar fields at January 2011 (month \# 1 for our trajectories). Highlighted regions are the locations of the finite-time coherent sets at January 2011.}
\label{fig:8rotvecs-hatfunc}
\end{figure}

The finite-time coherent sets are encoded in Figure \ref{fig:8evecs-hatfunc} as regions in the ocean that have extreme positive or negative values in one or more of the eigenfunctions $f^{(k)}$, $k=1,\ldots,8$.
Typically the spatial scale of the features highlighted in the eigenfunctions decreases as one proceeds down the eigenspectrum.
Some of these extreme value regions occur in more than one eigenfunction and so to produce a clearer set of images, where there is \emph{exactly one highlighted feature per function}, we apply the SEBA algorithm.
SEBA (Sparse EigenBasis Approximation) \cite{SEBA}\footnote{\textsc{Matlab} code is listed in \cite{SEBA}, and \textsc{Matlab} and {\bf \textsf{julia}} code may be downloaded at \url{http://www.maths.unsw.edu.au/~froyland/software.html}.} produces a sparse basis that approximately spans a given eigenbasis;  in the present case the span of $f^{(1)},\ldots, f^{(8)}$.
We input the eigenvectors $u^{(1)},\ldots,u^{(8)}$ to SEBA, and the algorithm outputs sparse vectors $s^{(1)},\ldots,s^{(8)}$ of the same size.
The sparse vectors $s^{(k)}$ mostly take values between 0 and 1 and the number $s^{(k)}_i$  can be interpreted as the likelihood that the trajectory $x_{i,t}, t=1,\ldots,T$ belongs to a coherent set.
From these sparse coefficient vectors, we may construct sparse functions $\hat{f}^{(k)}, k=1,\ldots,8$ as $\hat{f}^{(k)}=\sum_{i=1}^{I+C} s^{(k)}_i\varphi_i^1$, shown in Figure \ref{fig:8rotvecs-hatfunc}.


\section{Discussion of the results of the finite-time coherent set analysis}

\begin{figure}[hbt]
  \centering
  \includegraphics[width=\textwidth]{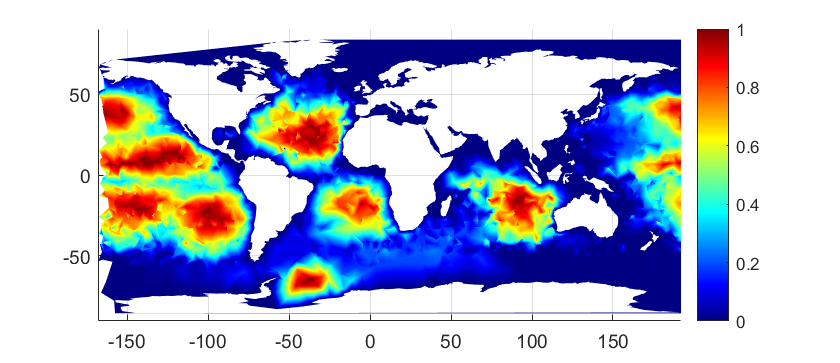}
\caption{Maximum function $f^{\rm max}$ formed from maximising the sparse basis functions $\hat{f}^{(1)},\ldots,\hat{f}^{(8)}$ shown in Figure \ref{fig:8rotvecs-hatfunc}, represented as a scalar field at January 2011 (month \# 1 for our trajectories). Highlighted regions are the locations of the finite-time coherent sets at January 2011.}
\label{fig:max8evecs-hatfunc}
\end{figure}

In Figure \ref{fig:max8evecs-hatfunc} we combine the eight coherent features displayed in Figure \ref{fig:8rotvecs-hatfunc} into a single image by maximising across the sparse vectors $s^{(k)}$.
Set $s^{\rm max}:=\max\{s^{(1)},\ldots,s^{(8)}\}$, where entry-wise maximisation is meant;  i.e.\ $s^{\rm max}_i=\max\{s^{(1)}_i,\ldots,s^{(8)}_i\}$ for each $i$.
We generate the corresponding maximum function as $\hat{f}^{\rm max}=\sum_{i=1}^{I+C} s^{\rm max}_i\varphi_i^1$.

\begin{figure}[hbt]
    \centering
    \includegraphics[width=\textwidth]{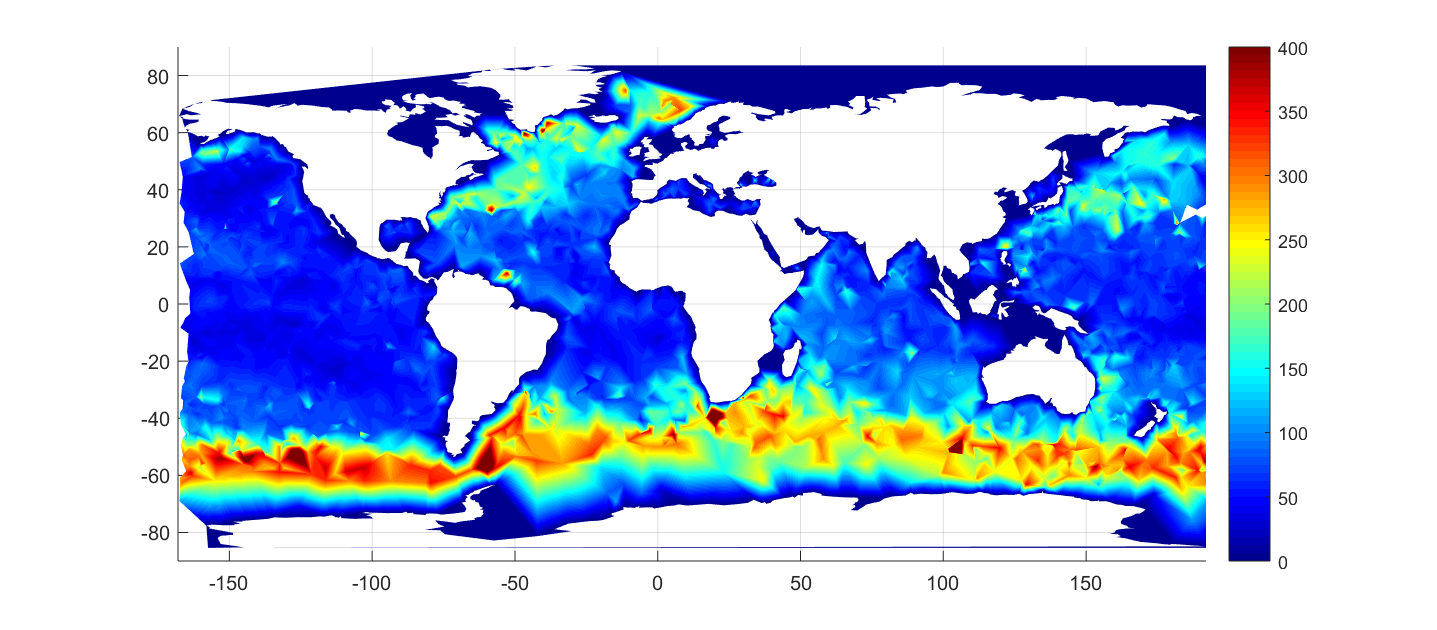}
    \caption{Plot of RMS speed in km/month, computed from Argo floats monthly across all 72 months, illustrated with basis functions $\varphi_i^{36}$ at month 36.}
    \label{fig:globalRMShatfunc}
\end{figure}

In Figure \ref{fig:max8evecs-hatfunc} we see as colour highlights the eight most coherent deep ocean objects over a six-year time frame;  these correspond to several known oceanographic features.
In the Atlantic, the 6-year coherent sets highlight the North and South Atlantic Subtropical Gyres in red, while the equatorial countercurrent and North Atlantic Current are excluded in dark blue.
In the North Pacific, the Subpolar Gyre is separated from the Subtropical Gyre, while in the South Pacific, the South Pacific Gyre is separated into eastern and western components.
In the Indian Ocean, the northern extent of the Indian Ocean Gyre is limited by stronger mixing due to the equatorial countercurrent.
In each of the Atlantic, Pacific, and Indian Oceans, the western extents of the highlighted features is limited by the stronger mixing of the western boundary currents.
The Weddell Gyre off the Antarctic coast is visible for the first 20 months of the 72-month flow, but after that time all floats in our dataset disappear from the Weddell Gyre and therefore its signature can no longer be detected.
In Figure \ref{fig:globalRMShatfunc} we show average velocity of the floats across the six-year period, computed directly from monthly spatial increments of the floats, and interpolated using the basis functions $\varphi^t_i$.
There is a partial negative correlation between average float speed according to Figure \ref{fig:globalRMShatfunc} and the likelihood of float membership in a coherent set shown in Figure \ref{fig:max8evecs-hatfunc}, but  fast currents need not imply the absence of coherent dynamics, nor need the converse implication be true.

\subsection{Estimating optimal boundaries of the six-year coherent sets}
The likelihood representation of Figure \ref{fig:max8evecs-hatfunc} is useful to gauge the relative coherence of regions, with red indicating greater coherence.
We can determine a boundary of each of these coherent sets highlighted in Figure \ref{fig:max8evecs-hatfunc} by using (\ref{cheeger}).
Recall that each $\hat{f}^{(k)}$ takes positive values on a coherent set.
We will use superlevel sets $A^k_c:=\{x\in X: \hat{f}^{(k)}\ge c\}$ as candidates for $A$ in (\ref{cheeger}).
Recall that we are interested in sets $A$ that minimise $h(A)$ because such regions have minimal average evolved boundary length relative to enclosed area, and therefore such a set will have minimal mixing in the presence of small diffusion.
By using level sets we reduce the minimisation problem to considering a one-parameter family of sets $A_c$ indexed by the scalar $0<c<1$ and minimise over $c$.
In practice, for a fixed $k$ we (i) discretise the range of $c$ into steps of 1/100, (ii) evaluate the right-hand-side of (\ref{cheeger}) for each discrete value of $c$, and (iii) select a local or global minimiser $c$, denoted $c_{\min}$.
In step (ii), for a given $c$ we must estimate the boundary length and area of the evolved $A^k_c$ at each discrete time $t=1,\ldots,T$, namely $\ell(\partial(\phi_t(A^k_c)))$, and  $a(\phi_t(A^k_c))$.

\textit{Estimating the length of the boundary of $\phi_t(A^k_c)$:}
To estimate $\ell(\partial(\phi_t(A^k_c)))$ for each $t=1,\ldots,T$ we interpolate the values  $s^{(k)}_i$ at the locations $x_{i,t}$ onto a fine, uniform spatial longitude/latitude grid (1 degree by 1 degree) and then use a contouring routine.
For example, with \textsc{Matlab}, we used \verb"contourcs" to obtain \verb"contourc" output as a struct array\footnote{\url{https://www.mathworks.com/matlabcentral/fileexchange/28447-contourcs-to-obtain-contourc-output-as-a-struct-array}}.
This contouring routine outputs a sequence of points so that the polygonal arcs connecting those points approximates the level set $\{x\in X: \sum_{i=1}^{I+C} s_i^{(k)}\varphi_i^t= c\}$ (recall the $\varphi_i^t$ are the piecewise affine hat function basis constructed from the mesh at time $t$), which is the boundary of $A^k_c(t)$.
The length of the boundary of the resulting polygon (or union of polygons) is then easily computed by the Euclidean lengths of its edges, scaled appropriately according to latitude.

\textit{Estimating the area of $\phi_t(A^k_c)$:}
To estimate the area of $A^k_c(t)$ we first apportion an area to each grid point in our 1 degree by 1 degree uniform grid, taking latitude into account.
We then simply sum those areas corresponding to grid points whose interpolated $s_i^{(k)}$ value is above $c$.

\begin{figure}
    \centering
    \includegraphics[width=0.8\textwidth]{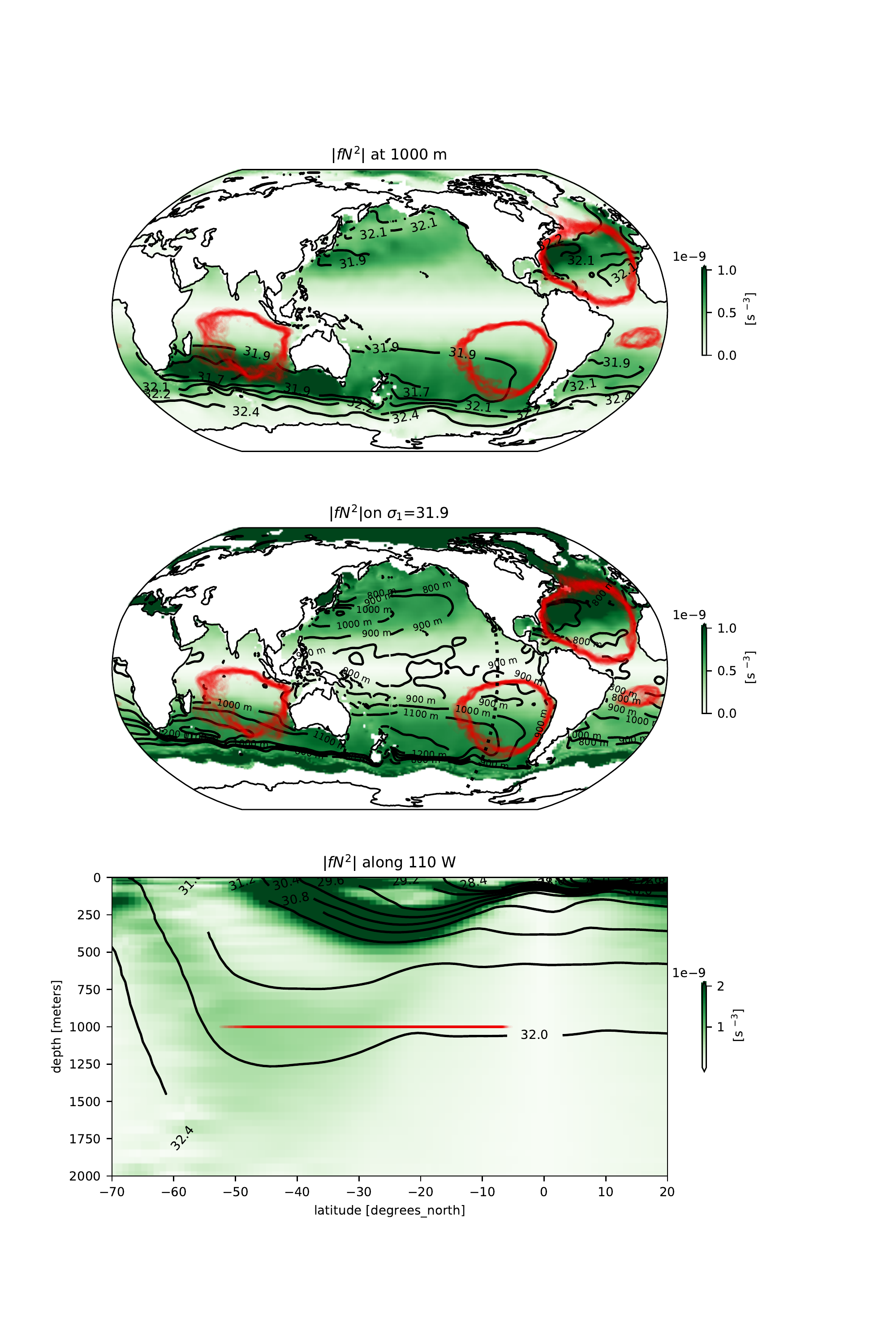}
    \caption{Boundaries of finite-time coherent sets (SEBA function contours) in the eastern South Pacific, North and South Atlantic, and Indian Oceans (red)  together with climatalogical absolute value of planetary potential vorticity $f N^2$ calculated from the World Ocean Atlas \cite{Locarnini2018_WOA_Temp,Zweng2019_WOA_Salinity}.
    The upper panel shows the potential vorticity at 1000 db (the Argo float parking depth), with contours of potential density $\sigma_1$.
    The middle panel shows the same field on the $\sigma_1 = 31.9$ potential density surface, with contours of depth.
    The lower panel shows a cross section through the middle of the SEBA coherent region in the eastern South Pacific, along the section indicated by the dashed line.
    Black contours are $\sigma_1$ isopycnals.
    }
    \label{fig:PV_spac}
\end{figure}

Putting the above boundary length and area calculations together, we have for a given threshold $c$, estimates of  $\ell(\partial(\phi_t(A^k_c)))$, and  $a(\phi_t(A^k_c))$ for each $t=1,\ldots,72$.
These values are input to (\ref{cheeger}), where the integral in (\ref{cheeger}) becomes a discrete sum across $t=1,\ldots,72$.
The result is a value $h(A^k_c)$.
We create these numbers for a discrete collection of $c$ between 0 and 1 and select a local or global minimiser $c_{\min}$.
Once such a $c_{\min}$ has been found, this defines our coherent set $A_{c_{\min}}$ at our initial time.
Having fixed $c=c_{\min}$, to construct the evolution of this coherent set over the 72-month time frame, we create a sequence of sets using the same contouring approach for each $t=1,\ldots,T$, using the \emph{same} threshold value $c_{\min}$.
These 72 boundaries for the coherent set associated with the eastern component of the South Pacific Subtropical Gyre, the North and South Atlantic, and the Indian oceans are plotted in red in Figure \ref{fig:PV_spac} (upper).

We emphasise that the concept of coherent sets is specifically designed to allow \emph{spatial motion} of the sets over time.
This is distinct from an earlier concept of almost-invariant sets \cite{dellnitz1999approximation}, where approximately \emph{fixed-in-space} objects, such as the Antarctic gyres, are mapped \cite{FPET07,dellnitz09}.
In the present experiment, we have relatively sparse Argo float information and our dynamic Laplacian approach allows us to reliably identify large basin-scale and sub-basin-scale coherent features.
While the concept of coherent sets allows these features to move in time, because of their large spatial scale in this application, they naturally tend to remain approximately in place, as illustrated by the red curves in Figure \ref{fig:PV_spac} (upper).


\section{Comparison with large-scale ocean tracers}

In this section, we look for relationships between the coherent sets detected by the dynamic Laplacian method and large-scale oceanographic tracers.
Our working hypothesis is that the coherent sets represent regions within which water parcels can exchange freely but are relatively isolated from exchange outside of the region.
A complicating factor, however, is that the Argo trajectories are isobaric and thus approximately constrained to the 1000 m depth plane, while water parcel trajectories are in general three dimensional.
Specifically, in the quasi-adiabatic ocean interior, we expect water parcel trajectories to lie on the neutral plane.

A central feature of the mid-depth ocean circulation is the existence of ``shadow zones'' on the Eastern boundaries.
Shadow zones are predicted by classic ventilated thermocline theory \cite{LuytenPedloskyStommel1983,Pedlosky1986} as regions for which the geostrophic streamfunction forms closed interior contours and does not connect with the surface, thereby isolating the shadow zones from surface ventilation.
One prediction of this theory is that, within such isolated regions, potential vorticity becomes homogenized due to eddy mixing \cite{RhinesYoung1982}.
If the coherent set detected by Argo floats corresponds with the shadow zone, we might then expect potential vorticity to be homogenized therein.

To test this, we examine planetary potential vorticity (PV), defined as $f N^2$, where $f$ is the Coriolis parameter and $N$ is the Brunt-V{\"a}is{\"a}l{\"a} frequency (proportional to stratification).
The frequency $N$ was calculated from climatalogical temperature and salinity fields provided by the World Ocean Atlas \cite{Locarnini2018_WOA_Temp,Zweng2019_WOA_Salinity}.
In Figure~\ref{fig:PV_spac} (top), we plot the absolute value of PV at 1000 m depth, together with the boundaries of the coherent sets in four major ocean basins (Indian, eastern South Pacific, North Atlantic, and South Atlantic).
The South Atlantic set is anomalously small, and therefore does represent a region of relatively homogeneous PV.
In the three largest sets (Indian, eastern South Pacific, and North Atlantic) rather than being homogenized within the coherent set, PV shows a strong meridional gradient across the region.
To test whether these gradients are an artifact of projection onto a depth surface, we also examine PV on an isopycnal surface that intersects the 1000 m surface in mid-latitudes ($\sigma_1 = 31.9$ kg m$^{-3}$); the same large scale gradients are present within the coherent sets.
This indicates that the dynamically isolated regions identified by the dynamic Laplacian analysis do not correspond closely with the classical shadow zones.
Moreover, textbook shadow zones---with homogenized PV along the eastern boundary of the subtropical gyre---are not evident in the PV field at these depths.

\begin{figure}
    \centering
    \includegraphics[width=0.8\textwidth]{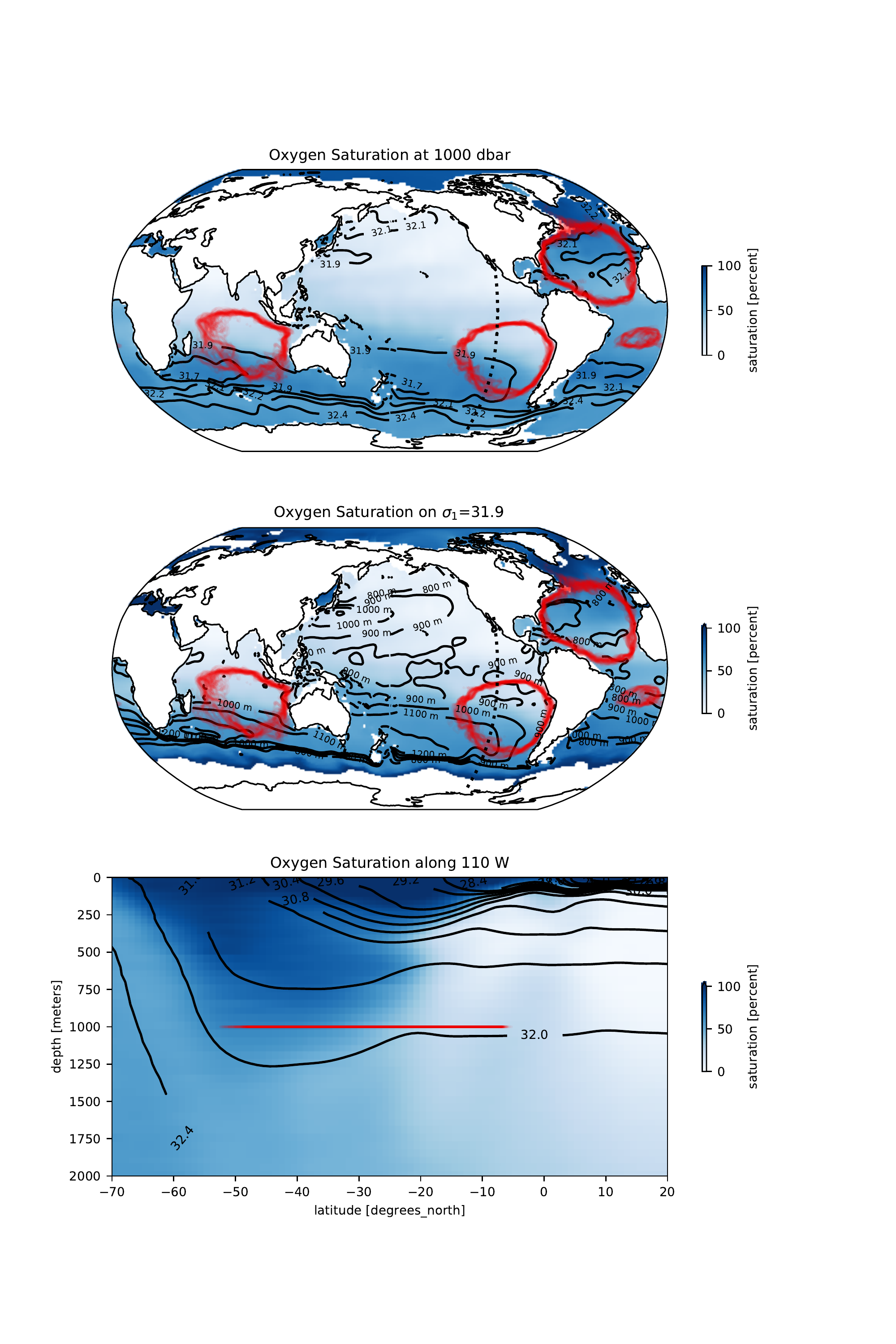}
    \caption{As for Fig.\ \ref{fig:PV_spac} except color shading now shows World Ocean Atlas Oxygen saturation \cite{Garcia2019_WOA_Oxygen} in blue.
    }
    \label{fig:O2_spac}
\end{figure}

Another prominent tracer feature at these depths are Oxygen Minimum Zones (OMZs).
OMZs are mid-depth layers of low O$_2$ concentration found in most ocean basins \cite{PaulmierRuizPino2009} and references therein.
OMZs are generally anchored to the Eastern boundary.
Their vertical extent begins near the surface but can reach as deep as 4000m in some ocean basins.
Most OMZs thus intersect with the 1000m depth surface sampled by Argo trajectories.
OMZs are formed through the combined effect of biology, circulation, and mixing; biological consumption is ultimately responsible for depleting oxygen, but the effects of circulation help determine the shape and position of the OMZs \cite{Wyrtki1962}.
In particular, regions of weak deep circulation and low mixing are associated with OMZs.
It is reasonable to ask, therefore, whether the dynamically isolated regions revealed by the dynamic Laplacian analysis correspond with OMZs.

To test this hypothesis, in Figure~\ref{fig:O2_spac} we examine the boundaries of the  coherent sets together with the oxygen saturation level from the gridded World Ocean Atlas dataset \cite{Garcia2019_WOA_Oxygen}.
Similarly to PV, there is no obvious overlap between low and / or homogenous oxygen saturation levels and the coherent sets.
Instead, the South Pacific and South Indian sets cross directly over the strong subtropical oxygen gradient, while the North Atlantic set occupies the entire basin.
The coherent sets do not correspond to the zonally elongated, Eastern intensified regions characteristic of OMZs \cite{PaulmierRuizPino2009}.

The overall lack of correspondence between large-scale ocean tracers (PV and oxygen saturation) and Argo coherent sets indicates that processes besides lateral transport at 1000 dbar must be playing a significant role in determining these tracer distributions.
Most notably the Argo floats experience no vertical transport, their ballasting confining them to the 1000 dbar surface.
Real water parcels, in contrast, move predominantly along neutral surfaces (isopycnals) at these depths.
While the local difference between the neutral angle and the 1000 dbar surface may be small at any particular point, over long trajectories, parcels originating at 1000 dbar may drift significantly away from this depth as they follow isopycnals over the basin scales representative of the coherent sets.
The potential mismatch between isobaric and isopycnal trajectories is therefore a major caveat around the coherent sets we have inferred from Argo trajectories.
Furthermore, diapycnal advection and mixing likely also play a strong role in tracer budgets at 1000 db, causing the observed tracer distributions to diverge from the patterns revealed by the Argo-based coherent sets.

\section{Conclusions}



The dynamic Laplacian approach \cite{F15} to identifying coherent sets and its FEM implementation \cite{FJ18} has previously been deployed on eddy detection in the Agulhas current \cite{FJ18} and in the North Atlantic \cite{SEBA}, where in both cases an altimetry-derived velocity field was used.
The eddy tracking experiment in \cite{FJ18} used trajectory sampling that was rich enough to compute reliable spatial derivatives.
This was relaxed in \cite{SEBA} below the level where reliable spatial derivatives of the flow map could be calculated, without impacting on the performance or resolution of the eddy tracking.
In the present paper we have dramatically reduced the input trajectory resolution to demonstrate that the dynamic Laplacian approach is well suited to the extremely spatially sparse trajectory data that is typically recorded from drifter datasets such as Argo.
Moreover, the approach easily accommodates the irregular temporal sampling of this dataset.

Using Argo trajectories over six years we identified the eight most coherent subsurface ocean features on this timescale (Figure \ref{fig:8rotvecs-hatfunc}).
Because of the global nature of the trajectory data, many of these dominant coherent regions are aligned with the major ocean boundaries, but with a bias toward the east of the basin due to the stronger mixing western boundary currents that destroy coherent water motion.
In particular, the vicinities of the Mozambique, Brazil, Kuroshio, and East Australian currents are conspicuously not highlighted in Figure \ref{fig:8rotvecs-hatfunc}.
There is a separation between the North and South Atlantic associated with the equatorial countercurrent and for the same reason the main coherent region in the Indian Ocean is pushed southwards.
In the North Pacific, the separation between the subpolar and subtropical gyres is evident, and the South Pacific gyre splits into eastern and western components.
To our knowledge, the dynamical separation between these two parts of the subtropical gyre has not previously been identified and merits further investigation.
In contrast to ocean eddies, which move considerably in space over time, the global-scale coherent regions we have identified move only slightly throughout the six-year duration, reflecting the relatively stationary nature of the ocean dynamics at these large spatial and temporal scales.

We hope that this analysis can help inform future deployment of Argo floats.
The coherent sets we have identified represent geographical regions from which floats are unlikely to escape.
Floats deployed near the center of the coherent sets are nearly certain to remain confined to the regions shown in Figure~6.
In contrast, floats deployed outside the sets, or near set boundaries, may be free to wander into different regions.
This information could help target float deployment to sample specific regions.

Another potential application of the method introduced in this paper would be to study the dynamic geography using other Lagrangian float ensembles.
The most obvious candidate would be the surface global drifters.
It would also be interesting to explore the regional dynamic geography of ocean basins with a high density of Lagrangian floats.
For example, the Consortium for Advanced Research on the Transport of Hydrocarbon in the Environment deployed over 300 drifters in the Gulf of Mexico in 2012 \cite{MarianoEtAl2016}, the so-called Grand Lagrangian Deployment, GLAD.
The high density of drifters in the Nordic Seas \cite{KoszalkaEtAl2011} makes this another attractive region to explore.

We hypothesized that the Argo coherent sets could be somehow aligned with large-scale ocean tracers and tested this by visualizing the relationship between set boundaries and isosurfaces of potential vorticity and oxygen saturation.
In general, both tracers and coherent sets show zonal asymmetry, with more coherence / tracer homogenization towards the Eastern boundary.
However, there was no clear smoking-gun relationship found in the tracers we examined.
Instead, the coherent sets mostly tended to occupy large regions, up to the entire basin scale, while the tracers contain notable features (e.g.~shadow zones, oxygen minimum zones) on smaller scales.
The coherent sets tended to overlie the strong subtropical gradient present in both tracers, suggesting that the sets do not represent regions of clear tracer homogenization.
The mismatch between the tracer isosurfaces and the coherent sets is likely due to three-dimensional effects (e.g. flow along isopycnals rather than isobars, overturning circulation) and non-conservative processes (e.g. mixing, oxygen consumption).
Future work examining coherent sets in three dimensions could help resolve these questions.
Such work would have to rely on numerical trajectories from an ocean model, as all existing real floats and drifters are constrained to 2D motion.




\section{Acknowledgements}
RPA acknowledges support from NASA award NNX 80NSSC19K1252.
GF thanks the Banff International Research Station for supporting a five-day workshop in January 2017 where initial discussions for this work took place.
The research of GF has been partially supported by two ARC Discovery Projects over the course of this research.
KS was supported by an ARC Discovery Project.

\newcommand{\etalchar}[1]{$^{#1}$}

\end{document}